\documentclass[12pt]{article}
\usepackage{fullpage, amsfonts}

\author {A.~M.~Vershik\thanks{%
St.~Petersburg Department of Steklov Institute of Mathematics.
E-mail: vershik\@pdmi.ras.ru.}}

\date{June 2007}

\title{L.~V.~Kantorovich and linear programming}

\begin{document}
\maketitle

I want to write about what I know and remember about the activities
of Leonid Vital'evich Kantorovich, an outstanding scientist of the 20th
century; about his struggle for recognition of his 
mathematical economic theories; about the initial stage of the history
of linear programming; about the creation of a new area of mathematical activity
related to economic applications, which is called sometimes operation research,
sometimes mathematical economics, sometimes economic cybernetics, etc.; about its
place in the modern mathematical landscape; and, finally,
about several personal impressions of this distinguished scientist.
My notes in no way pretend to exhaust these topics.

\section{The ``discovery'' of linear programming}

When attending a remarkable detailed two-year course of
functional analysis taught by L.~V.~Kantorovich (in 1954--1955), I had
never heard neither about his research in duality theory,
nor about computations of Banach norms (his notes in {\it Doklady
Akad.\ Nauk} published in 1938--1939), nor about linear extremal
problems (the famous ``veneer trust'' problem) and the method
of resolving multipliers suggested by him for solving problems that
were later called problems of linear programming. I learned all
this a little later. As to this course of functional analysis,
L.V.\ taught it for several years; later it formed
the basis of the widely-known book 
{\it Functional Analysis in Normed Spaces},
written by L.V.\ and  G.~P.~Akilov, his
main pupil in this field. At the time, it was undoubtedly one of the most
exhaustive and deep monographs  on functional
analysis
in the world literature, and simultaneously a textbook. 
Later I had a chance to
make sure that abroad it was also very popular. By the way, the ``Leningrad'' functional
analysis, initiated by V.~I.~Smirnov, G.~M.~Fikhtengolts,
L.V., who served as the main engine, and, somewhat later, G.~P.~Akilov, had its
specific feature: the influence of mathematical physics (S.~L.~Sobolev),
complex analysis (V.~I.~Smirnov), theory of functions (G.~M.~Fikhtengolts,
I.~P.~Natanson, S.~M.~Lozinsky) was stronger than in the Moscow
or Ukrainian schools, which were more affected by operator theory, spectral theory, multiplicative
functional analysis, representation theory, and Banach geometry.
Before the war, L.V.\ also created a specific ``Leningrad'' research direction:
functional analysis in partially ordered spaces. But the main contribution
of L.V.\ in this field, unanimously acknowledged throughout the world,
is related to applications of functional analysis to approximative methods
(summarized in his famous paper ``Functional analysis and approximative
methods'' published in 
{\it Uspekhi}\footnote{The journal {\it Uspekhi Matematicheskikh Nauk}.}). These works
were recognized by the Stalin Prize; they initiated an enormous amount of
research in this area.

During many post-war years, problems
of functional analysis were mainly being discussed at the well-known 
Fikhtengolts--Kantorovich seminar at the Department of Mathematics
and Mechanics of the Leningrad State University, which I attended
regularly since 1954 and up to its actual closing in the mid 50s.
An active role in its organization, 
especially in the last years,
was played by Gleb Pavlovich Akilov, my first scientific
advisor, an original and independent man, 
a pupil, coauthor, and colleague of L.V. Once 
G.~Sh.~Rubinshtein, who was also in fact a pupil of L.V., 
gave a talk 
on best approximation and the problem concerning the structure of the
intersection of a ray with 
a cone, i.e., essentially a problem of linear programming. But at the time
it was perceived merely as a separate talk on a particular subject, 
and I do not remember L.V., or somebody else, giving any comments or
saying in what context it should be understood. But I remember the
impression of reticence remained after this talk.

Apparently, this reticence was due to the internal veto, whose
reasons were well known to senior participants of the seminar,
implicitly imposed on open discussions of this circle of L.V.'s papers.
This veto ensued from the persecution of his ideas
unleashed by the ideological bonzes soon after he had published
the brilliant
booklet {\it Mathematical Methods
in the Organization and Planning of Production} (1939)
and written, during the war,
a book on economics, which was published almost 20 years later.
This turn of events threatened 
to bury the whole research direction, as well as 
to bury his author, in the most direct sense.
Only many years later it became known how serious were the accusations
and threats of high scientific and ideological officials. This veto 
existed up to 1956. And it applied not only to economic matters,
but even partly  to the mathematical
aspect of L.V.'s works. Many of these documents have been recently 
discovered by V.~L.~Kantorovich. It is of great importance to bring them
into the open for all who are interested in the history
of our science. Already at that time there were some vague 
conversations about the applied research carried out by L.V.\ in 
the post-war years:  on optimal
cutting with V.~A.~Zalgaller,  on the mass transportation problem with M.~K.~Gavurin, etc.
But, honestly speaking, 
I attributed all this to the boring category of ``collaboration of science and 
industry'' (a propaganda clich\'e of the time, which usually
covered superficial or even meaningless works) and did not know
the economic and mathematical seriousness of the matter.
During the first years, V.~A.~Zalgaller, M.~K.~Gavurin, and
G.~Sh.~Rubinshtein (one should also include in this list A.~I.~Yudin,
a student of L.V.\ who died
at the front, and maybe some others) were the closest
L.V.'s assistants in applied economic studies and developed the 
theory of these problems. With M.~K.~Gavurin, L.V.\ wrote
the famous paper on the mass transportation problem  (written before the war,
but published only in 1949).
With V.~A.~Zalgaller, he worked on the problem of optimal cutting and
wrote a book on this subject (1951); V.A.\ also tried to introduce 
optimal cutting at the Egorov
Wagon Plant in Leningrad. For well-known reasons, people
with ``flawed'' biographies could find a job at
civil enterprises (such as the Egorov Plant).  Sometimes
this resulted in the professional level at such a plant
being above the average. For the same reasons, G.Sh.\ (patronized by L.V.)
managed to get a job at the Kirov Plant, where he also tried to
introduce optimization methods and simply reasonable approaches 
to problems of local planning. Note that G.Sh.\
graduated from the university at the time when he, being a veteran of 
the war and a successful student, could not enter
a graduate school. Before the war, G.Sh.\ studied 
at the Odessa University and was a student of M.~G.~Krein, so that he
successfully combined the knowledge of that part of the work of M.~G.~Krein
and the whole Ukrainian school of functional analysis 
that was close to L.V.'s activities  (the L-moment
problem) with the good understanding
of L.V.'s ideas in linear programming. There were also attempts to
introduce optimization methods at the Skorokhod Factory, 
Lianozovo (ex-Egorov) Wagon Plant, Kolomna Locomotive Plant, etc.
But, amazingly, this activity met the resistance of those who would seem
to be most interested in it. At the time, as well as later, there existed
a number of comical examples of reasons why some or another well-founded
suggestion did not find support. For instance, suggestions on optimal
cutting came into conflict with the bonus promised to those who
collected more waste products for recycling, etc. Afterwards,
optimal cutting was much studied by the Novosibirsk pupils of L.V.,
in particular, E.~A.~Mukhacheva.

Were there any serious reasons why this useful activity met such
difficulties and was not on demand at the time? All of the few
papers on this subject written at those ``underground'' years 
were meant for engineers and published in nonmathematical editions
available for engineers. It would seem to be the best example
of ``cooperation of science and industry,'' which opened wide horizons
for scientific, mathematically-based, local and global
planning in economics. At the early period (1939--1949), one might think
that the reason of this antagonism was 
in the unpreparedness of people and their working conditions 
to comprehend these ideas and techniques, as well as in the deathening
ideological dogmata and stupidity of the Party supervisors and ideologists.
One might think that if the authorities were more enlightened,
they would be able to appreciate, implement, and use the new ideas.
Perhaps, L.V.\ also thought so. But the entire subsequent Soviet
history showed that the situation was much worse{\ldots} At the time, and even
later, it was not thoroughly understood that the reason behind the failure of
implementing the most part of new economic (and other) ideas
was not in particular circumstances or in the stupidity of bureaucrats,
but in the fact that the whole Soviet economic system, or, as it was
called later, command-administrative system, is organically unable
to accept any innovations, and no serious, big or small, economic reforms
leading to stability just cannot be realized inside it. 

It was not until the middle of 1956 that L.V.\ started actively promoting
this subject and giving lectures at the Department of Mathematics and Mechanics
and other departments of the Leningrad University, at 
the Leningrad Department of the Steklov Institute of Mathematics (LOMI). This was the
disclosure of a new, earlier forbidden, area. He told about the contents
of his 1939 book, about resolving multipliers, various models and
problems, etc. For the overwhelming majority of the audience, including myself,
these topics were completely, or almost completely, new. Undoubtedly,
the ``declassification'' of this subject had to do with the new hopes appeared
after Stalin's death, Khrushchev's speech, and the beginning of the ``Thaw.''
At this point, it is pertinent to recall V.~I.~Arnold's story about A.~N.~Kolmogorov:
being asked by V.I.\ why in 1953--1954 A.N.\ suddenly started working
on the classical and very difficult problem of small denominators
(this was the beginning of what is now called KAM theory), on which
he had never worked before, A.N.\ answered: ``Some hope has emerged.''

Undoubtedly, some hope had emerged also for L.V., the hope that he would finally
be able to explain and implement his ideas and to
overcome the Soviet dogmatism and obscurantism in economics.

When one says that the Soviet science (not the whole science,
but, say, mathematics) was successfully developing and reached
a very high level, there is no arguing, but one still must not forget about
this and many other similar stories:
ideological pressure, selection according to personal details,
etc.\ never allowed a gifted man to realize his talent
in full measure or even did not allow to realize it at all.
The indubitable scientific 
achievements of the Soviet period is only a small part of what could 
appear in the conditions of freedom, and the losses caused by
frustrated or forbidden discoveries and ideas are irreparable.

At this period (the late 50s --- the early 60s) L.V.\ developed an enormous
activity. His numerous vehemently delivered lectures, polemical talent and
zeal ignited the audience. I remember the intellectual attack organized
by L.V.\ (in about 1959) in connection with taxi tariffs. 
He was charged with this task by some
authorities (apparently, this was a test). L.V.\ organized a team consisting of
15-20 mathematicians, each being given a separate problem. The 
situation required a brainstorm: the team had to 
carry out a detailed analysis of a pile of data
and output its recommendations within a week.
There were some exaggerations --- sometimes L.V.\ could
be carried away by his ideas
and suggest unrealistic projects, but the task was
fulfilled and L.V.'s recommendations on taxi tariffs (for instance, the
idea of initial fare) were implemented in 1961 and have been used since then; moreover,
L.V.'s predictions (the results of investigating the elasticity of demand)
proved to be absolutely accurate.

Mathematicians listened to L.V.'s talks and series of lectures with
enthusiasm. Gradually, the number of those who mastered these
techniques at LOMI and at the Department of Mathematics increased. At first, 
the then Head of the Department S.~V.~Vallander took an 
active part in the popularization of L.V.'s ideas.
The Department organized
a series of lectures by L.V.\ for a wide audience. He also gave many
talks at the LOMI General Seminar.

But L.V.'s lectures for economic audience encountered hostility
or, in any case, scepticism. I remember comical and ignorant
objections of political economists during his lectures at the Department of Economics.
After the well-known Khrushchev's speech, the ideological blinders were
somewhat weakened, and it became more difficult to defend stereotyped
rubbish. One could see that the orthodox position was weakening, and
among  political economists and ideologists, people began to appear who
wanted to understand. Once (in 1957) I encountered, 
in an unofficial situation,
G.~V.~Efimov, an orientalist who was at the time the Vice-Rector for Research of
the Leningrad University, a man
not of a liberal type. However, to my surprise, he was captivated by
my description of L.V.'s ideas and their possibilities as they seemed
at the time.

The idea that 
turned out to be the most important for the whole economic
theory is the direct economic interpretation of the dual problems formulated
by L.V., and it is this idea that provoked the hostility of orthodox economists.
The economic analog of the variables of the dual problem
(resolving multipliers), which afterwards L.V.\  aptly called ``objectively
determined valuations,''\footnote{In the English literature, they are
known as ``shadow prices.''} was, roughly speaking, the exact mathematical equivalent
of the notion of prices, and one should call them so if it were not for
the fear of ideological invectives of that time. The subtlety of the term 
invented by L.V.\ (``objectively
determined valuations'') was in the fact that, funny as it may seem, 
Marxists are defenseless against the term ``objectively.''
The emphasis placed by L.V.\ on dual problems led to the most important
economic conclusions and protected the common sense against
standard dogmata, in particular, advocated natural rent, true
estimate of costs, etc. This was his most important contribution, 
his main trump card in discussions, and the most irritating issue
for his opponents, who, of course, accused him of revising Marx's
``labor'' theory of value, even more so because labor
appeared in his theory, too, and had no difference from, say, any raw material.
How much effort he exerted to defend himself against these fatuous
attacks! One could write a book about this using the documents stored in his
archive. Even A.~D.~Alexandrov, the then Rector of the
Leningrad University, could not (out of prudence or following direct
instructions) publish L.V.'s new book on economic calculation
at the University Publishing House.

Here is another example of how the officials of that time were afraid
of everything related to this subject. At about the same time (1957), 
I, with a coauthor, wrote a popular paper on mathematical economics
for the newspaper {\it Leningradskaya Pravda}, having a preliminary
agreement with a member of the editorial board, who was an acquaintance of mine.
Nevertheless, the paper was never published. Having smelt something
nonstandard, the editor required that the text, which
was no more than a popular paper, should be agreed with the authorities,
which I refused to do.

One can judge
to what extent L.V.'s works were known to the scientific community
by the following fact. Once, at the end of 1956, G.~Sh.~Rubinshtein
wrote for me on a small piece of paper (I still keep it) 
the list of {\bf all} literature
in Russian on this subject, and this list contained only 
five or six items, including
L.V.'s 1939 booklet and his joint book with V.~A.~Zalgaller on
optimal cutting. Moreover, almost all of them were published in
little-known and rare editions, and nothing (except two or three L.V.'s
notes in {\it Doklady}) appeared in mathematical journals. It is interesting that
in the well-known volume {\it Forty Years of Mathematics in the USSR} (1959),
where the corresponding section was written by L.V.\ together with 
M.~K.~Gavurin, this topic occupies only one page and the same five
references are given. In spite of all this, these were years of hope
that progress, positive changes, and nondogmatic attitude towards
new ideas were possible.

As often happened in the USSR, military specialists were the
first to have an access to books 
on linear programming (Vajda), operation research (Campbell), etc.,
obtained
through special channels, translated into Russian, 
and not yet published in the USSR.
The interest of the military to this circle of problems
had to do not with economic problems (e.g., distribution of resources),
though they were also of importance for them; but
they were a part of the general theory of control of systems, which
was afterwards called by the strange term ``operation research.'' Undoubtedly,
in those years many scientific ideas obtained additional support if
for some reason they interested the military; operation research and,
in particular, linear programming, is one of such examples.

Of course, none of the military specialists (among which 
there were engineers that knew
mathematics quite well; some of them joined the Army after
graduating from mathematical and physical departments of universities)
had ever heard of L.V.'s works, and this is not surprising.
I remember how,  at the beginning of 1957, being on a business trip to the 5th Research Institute
of the Ministry of Defense at Moscow, I told to
D.~B.~Yudin and E.~G.~Golshtein, mathematicians working at this institute,
about resolving multipliers and L.V.'s works, and showed them the above-mentioned 
small list of references. 
For them, who had just started to
study the American literature on linear programming,
this was a revelation.
Later, they became the main writers on this subject, and their role
in its popularization is quite important. Indirectly,
their activity became possible exactly because they had  once been involved
in military research.

In the autumn of 1957, I asked L.V.\ to give a lecture at the 
Navy Computer Center, where I worked at the time. This large
center was created in 1956 along with another two ones, the Land Forces
Computer Center in Moscow
and the Air Forces Computer Center
in Noginsk near Moscow, on the tide
of the rehabilitation of cybernetics and the belated understanding of the
necessity to introduce first computers and modern 
mathematical and cybernetic methods into the Army. At this center there were
many serious specialists in automatic control, theory of shooting, and
other fields of military science. L.V.\ gave a 
public lecture on solution of some extremal problems,
which was a success. One of its
consequences was that military specialists, who had been thus far using
only foreign sources obtained through their own channels, 
started to believe that in this field too the works of our mathematicians
were pioneering. It was curious to see another piece of evidence that, in spite of 
the decades-long brainwashing about the priority of Russian and
Soviet science (and, most likely, just because of this
brainwashing), most people,
for example, many military specialists I met, on the contrary, could not
believe that something might have appeared in the USSR earlier than in the 
West. The humour of the situation was, exactly, that we had exchanged
our roles: they, being well-trained in ideology, in every lecture
harped on the ridiculous nonsense about priority.
However, they were highly sceptical when I told them about the
undoubted priority of L.V. And their sceptical attitude was quite understandable:
they hardly believed in the common phrases about Russian and Soviet
priority.\footnote{At this point, one cannot help remembering the sad story of I.~M.~Milin, a well-known
mathematician who worked at a military school in Leningrad and was
fired only because at a lecture, after the obligatory mention
of the priority of Russian mathematics in some elementary question,
he took the liberty of making a humorous remark: ``And now let's turn to
serious matters.''}

On the other hand, everybody knew quite well that
most often new and sensible ideas
appearing in the USSR could not push through, or managed to
push through only after a round-the-world voyage. This was partly the case
with L.V.'s theory, as well as with many other ideas.

L.V.'s offensive, which had begun in 1956, continued up to the mid 60s,
when his economic theories if not started to 
be recognized, but at least stopped to be forbidden
by the ideological and economic establishment.

Later, they even found some (not absolute) recognition: in 1965, L.V.\ 
(together with V.~V.~Novozhilov and V.~S.~Nemchinov) was awarded the Lenin Prize.
From the very beginning, L.V.\ obtained support ---
in discussions,
at conferences, etc.\ --- from many venerable
mathematicians (A.~N.~Kolmogorov, S.~L.~Sobolev) and some economists.
Many specialists participated in these discussions,
and, of course, they involved not only L.V.'s theories, but also many
other topics (close economic theories, e.g., that of V.~V.~Novozhilov,
cybernetics, the role of mathematics and computers, etc.). I remember
a large conference of mathematicians and economists 
that took place in 1960 in Moscow, where
many scientists, both venerable and young, spoke, and,
with rare exceptions, in support of the new ideas. As a whole, it
was undoubtedly the victory of sense, but L.V.\ wasted too much effort,
taken away from mathematics and science as a whole, on this struggle. 
In fact, from the late 50s, when
one of
his latest mathematical papers was published in {\it Uspekhi},
he stopped working systematically in ``pure'' mathematics.

The history of L.V.'s struggle for recognition of his ideas is extensive and
interesting both for a historian of science and a historian of the Soviet
Union. It is poorly reflected in the literature, and, unfortunately, 
few people study it now. Meanwhile, both this experience and
the economic principles advocated by L.V.\ are of importance today.
It is not until 1998
that the volume {\it Essays in the History
of Informatics in Russia} [in Russian], which, in particular, contains information
on this epic,  was published by the Siberian Branch of the Academy of Sciences
(Novosibirsk).

In 1989, we organized a conference in Leningrad dedicated to the 50th
anniversary of his classical booklet {\it Mathematical Methods
in the Organization and Planning of Production}. An account of this
conference was published in the journal
{\it Ekonomiko-Matematicheskie Metody}.
When preparing for this conference, V.~L.~Kantorovich found
many interesting and unknown documents related to L.V.'s struggle for
his ideas and, in particular, letters and instructions of ideological bonzes
concerning his works. These documents must be published and become
known to all who are interested in the sad and instructive history
of our country. At that time, and especially now, all these facts 
are known very poorly.

Of course, the awarding of the Nobel Prize put L.V.\ in a unique position
in the USSR (this was the only our Nobel Prize in economics, and,
moreover, awarded simultaneously
with the Nobel Peace Prize to A.~D.~Sakharov). Should it not mean 
general recognition
and confidence? However, up to the end his position remained
that of a captive rather than a principal expert, as it should be.

Though L.V.'s economic ideas were in a sense consonant with planned
economy and could be easily interpreted in generalized Marxist
spirit, the persistent hostility towards them
should be explained in psychological rather than
logical terms: the ignorance inherent in an ageing dogmatic regime
is not psychologically prepared for intellectual renovation, no matter how
clearly one tries to explain its advantages for the regime itself. A very simplified
interpretation of the mutual relations between L.V.\ and the ruling
ideology was given by A.~Katsenelinboigen in a quite interesting paper
``Does the USSR need Don Quixotes?'' ({\it L.~V.~Kantorovich:
a Man and a Scientist; His Contradictions} [in Russian], Chalidze Publications, 1990). 
I am not going to discuss deep and important problems related to the
mutual relations between a scientist and society. In the  Soviet period
these relations were especially complicated and did not allow for straightforward
and primitive interpretations. Of course, every conformist society
rejects new, unusually looking, ideas unless they are inculcated by
the authorities. This applies even if the advantage 
of the realization of these ideas is obvious. As a
French Sovietologist said on a similar occasion,
``The authorities
do not like to be protected by means they do not comprehend.'' No wonder that
a scientist wishing to advance his ideas is compelled to speak, at
least partly, in conformist language. And in this respect L.V.\ sometimes went too far.
But only those knowing and remembering that time and those people,
who had survived the chilling horror of the late 30s, can correctly
estimate some steps which would look strange in a normal human society. 
One should not dismiss the atmosphere of danger to life for those
who dared to deviate a little from the prescribed ideological dogmata,
and it is in this atmosphere that the generation of L.V.\
had spent the most part of life.

The famous paper ``Marx, Kantorovich, Novozhilov'' published in
{\it Slavic Review} by \break
R.~W.~Campbell demonstrated that some American economists
had a quite good understanding of what was going on in the USSR
with L.V.'s and V.~V.~Novozhilov's theories. This paper made quite a stir.
It was classified secret and kept in special depositories of
public libraries. And the authors (in particular, L.V.) had to
prove that they did not agree with the ``bourgeois'' interpretation
of theories and events suggested by Campbell. But in fact he 
gave a quite accurate description of the futility of the 
Soviet economic establishment, as well as the logical inevitability of the conclusions
that L.V.\ had derived by systematically developing his {\it strictly mathematical}
approach to concrete economic problems.

In the 90s, on many occasions 
I had to retell the epic of linear
programming in the USSR abroad, and it was extremely difficult to explain,
even by this example, the ``wonders'' of the Soviet system, which
rejected the achievements of its own scientists because of absurd
ideological prejudices. Perhaps, only the reference to the story of T.~D.~Lysenko,
well known in the West, helped the audience to understand at least something.

I would like to make another general remark. When we recall the
biographies of really outstanding Soviet scientists, we are
threatened by two extremes. The first one is to turn them into an icon, to 
remember only their scientific achievements and good deeds and forget
about their compromises with the authorities (such as signing obsequious  letters,
participating in ``collective'' campaigns, etc.). The second extreme is
to accuse them of open subservience to the totalitarian regime
because of the very usefulness of their work for the society. 
Now, when one can write openly,
when there is no censorial pressure, it is especially important
to understand that for many (not all) outstanding scientists of that
generation, their position in the Soviet society was, if not an inner
tragedy, at least a source of agony. Therefore no one of these
extremes allows one to fully understand the very complex
and tragic nature of this situation, the position of a
talented man under the pressure of total control. Some deeds
of these people may be regretted,
but it is not merely that their scientific achievements 
outweigh all the other things; one must also remember that the life of
a talented Soviet scientist was, above all, devoted to science, 
and sometimes, for the sake of science and the realization of his ideas, he had 
to compromise with the regime, which used his authority 
for its momentary
purposes and usually did not understand even its own advantage
from his activity as a whole, treating him, unless he had entirely become
its property or its adherent, with suspicion or even
hostility.

Turning back to linear programming, I think that the story of how the veneer trust
problem considered by L.V.\ in 1938 led to the theory of optimal
resource allocation is one of the most remarkable and instructive stories
in the history of the science of the 20th century. The same story can serve
as an apology for mathematics. This attitude towards L.V.'s works
has gradually become common among mathematicians; it was shared by
A.~N.~Kolmogorov, I.~M.~Gelfand, V.~I.~Arnold, S.~P.~Novikov, and others.
One cannot help admiring the naturalness and the inner harmony of
L.V.'s mathematical works on duality of linear programming
and their economic interpretations.

\section{Mathematical economics as a branch of mathematics and 
its connections with other fields}

\subsection{Connections of linear programming with functional and convex
analysis}

Already before the war, L.V.\ was a recognized authority in many fields
of mathematics, especially as one of the founders of a school of
functional analysis. No wonder that linear
programming  in his interpretation was also related to functional analysis. The same viewpoint
on these problems was shared by J.~von Neumann: his fundamental theorem of
game theory, models of economics and economic behavior, and other
results in mathematical economics have a strong flavor of conceptions
of functional analysis and duality.

Like most part of the mathematicians belonging to L.V.'s school, I
originally conceived the mathematical aspect of 
optimization econometrics
in the functional-analytic framework.
In other words, the duality scheme was being considered, in a natural way, in terms of
functional analysis. Undoubtedly, there is nothing more acceptable
from the conceptual viewpoint. Convex analysis, which was developed
after the 50s on the basis of optimization problems, gradually
absorbed a significant part of linear functional analysis, as well as
of the classical results of convex geometry. It is in this way that
I constructed 
my course of extremal problems that I taught at the Leningrad
State University for 20 years (1973--1992): it included general
(infinite-dimensional) separability theorems,
the duality theory of
linear spaces, etc.

Historically, the first known connections of L.V.'s theory were connections
with the theory of best approximation and, in particular, with M.~G.~Krein's
work on the L-moment problem. M.~G.~Krein was one of the first to
notice this. The real consequence of this observation was
the gradual understanding that
the methods for solving both problems are essentially similar.
The first such method  goes back to Fourier.
Later, in the 30s--40s, important studies were carried out by T.~Motzkin
and M.~G.~Krein's Ukrainian school (in particular, S.~I.~Zukhovitsky and
E.~Ya.~Remez). However, the method of resolving multipliers
and the simplex method were new for the theory of best approximation.
From the point of view of principle, the crucial issue was
the very interpretation of the problem of Chebyshev approximation as
a semi-infinite-dimensional problem of linear programming. Infinite-dimensional
programming was also the subject of several research works of my students at
the Department of Mathematics and Mechanics of the Leningrad State University
(M.~M.~Rubinov, W.~Temelt) and some Moscow mathematicians (E.~G.~Golshtein
and others).

The duality theory of linear spaces with cones provides a natural
language for problems of linear programming in spaces of arbitrary
dimension. Paradoxically, this was grasped by N.~Bourbaki, who was far from
any applications. Looking closer, in the exercises of 
the 5th volume of {\it Elements of Mathematics}
(rather an abstract opus!) one can find
even the alternative theorem for linear inequalities and a number
of facts close to duality theorems of linear programming. And this
is natural. Such fundamental theorems of classical linear functional analysis as
the Hahn--Banach theorem and linear separation theorems
are purest convex geometric analysis. The same is true for general
duality theory of linear spaces.

The classical Minkowski--Weyl theory of linear inequalities appeared in modern form
in H.~Weyl's work of the 30s, a little bit earlier than L.V.'s theories,
and this link is especially transparent. Alternative theorems,
Farkas lemmas, the Fenchel--Young duality in the theory of convex
functions and sets, all this was not combined with the theory of linear
programming until the 50s. 
Apparently, L.V.\ did not learn about all these links till somewhat later, 
but his contribution
is that he found a unified approach, based on the ideas of functional
analysis and revealing the heart of the problems. Simultaneously, this
approach provided a basis for numerical methods of solving these problems. One can say without
exaggeration that functional analysis became the basis
for the whole mathematical economics. A great number of problems of 
convex geometry and analysis (from Lyapunov's theorem on the convexity
of the image up to the theorems on the convexity of the moment map) 
are also related to these ideas
and their generalizations.

Later there appeared related works on the theory of linear
inequalities (S.~N.~Chernikov,
Ku Fan),
convex geometry, etc., whose
authors did not always know about the previous results; and now one still
cannot say that the research  in this area  is properly summarized. 

\subsection{Linear programming and discrete mathematics}

However, linear programming has serious relations with discrete mathematics
and combinatorics. More exactly, some problems of linear programming
are linearizations of combinatorial problems. Examples are the
assignment problem and the Birkhoff--von Neumann theorem, the Ford--Fulkerson
theorem, etc. At first, this aspect of the theory was not noticed by our
mathematicians, and its understanding came to us from the foreign literature 
somewhat later.
The fundamental theorem of the theory of zero-sum matrix games 
(the minimax theorem) was brilliantly connected to linear
programming by J.~von Neumann, see G.~Dantzig's memoirs cited in the paper
``John von Neumann'' by A.~M.~Vershik, A.~N.~Kolmogorov, and Ya.~G.~Sinai 
in vol.~1 of the book J.~von Neumann, {\it Selected Works on Functional
Analysis} (Nauka, Moscow, 1987). Dantzig writes about his 
striking conversation with von Neumann. Within an hour, von Neumann
described the connection between duality theory and theorems on matrix games
and outlined a method for solving problems
of linear programming. This connection
was not mastered at once. I remember that at first, the Leningrad specialists
in game theory
did not take into account that 
the problem of finding the solution of a zero-sum matrix game
is a problem of linear programming; the (undoubtedly beautiful) method for
solution of games belonging to J.~Robinson was considered almost the
only numerical method for finding the value of a game. The final proof
of von Neumann's minimax theorem (the first proof was topological
and used Brauer's theorem) essentially contained duality theory.
Later, the equivalence of a game problem and a problem of linear programming
has been intensively exploited.

In the first years, emphasis 
on the connections with discrete mathematics and combinatorics
prevailed in the most part of foreign articles on linear programming, 
while Soviet authors emphasized
the connection with functional and convex analysis and developed
numerical methods.

From the point of view of linear and convex programming, 
the most important parts
of combinatorics are the combinatorial geometry of convex
and integer polytopes and the combinatorics of the symmetric group.
The most important works of the first period
include a book by B.~Gr\"unbaum and papers by V.~Klee
on the combinatorics of polytopes,
and the studies of G.-C.~Rota and R.~Stanley  in general and algebraic combinatorics.
Simultaneously,
close topics appeared in singularity theory (Newton polytopes),
algebraic geometry (toric varieties and integer polytopes), etc.
Later, extensive connections were discovered with the symmetric group and
the combinatorial theory of Young diagrams, one of the main subjects
of the ``new combinatorics,'' as well as with posets and matroids.
It is interesting that I.~M.~Gelfand, who called combinatorics the
mathematics of the 21st century,
almost simultaneously (and independently)
arrived at a number of close problems (matroids,
Schubert cells, secondary polytopes). At present, new combinatorial problems 
play a key role in various fields of mathematics. 

In the first years, my interest to linear programming arose quite
independently from my mathematical favorites of the time and, in
particular, not only because I was a student of L.V.\ in functional
analysis and listened to his first exciting talks on linear programming
and its applications in economics. At that moment, this interest was rather practical
than theoretical. The point is that, having graduated from
the university and refused, for some reasons, to enter a postgraduate
school, I worked at the Navy Computer Center and got interested in the
problem of multi-dimensional best approximation as an applied researcher.
One of the problems I worked on
at the Computer Center was computer representation of firing
tables, and I suggested to approximate them instead of 
storying in the memory of a computer. I formulated a kind of generalization
of the problem of best approximation, namely, 
the problem of piecewise polynomial
best approximation (we had then never heard about splines) for functions
of several variables. Later, when I started to work at the university,
my first graduate students investigated this problem. Even more later,
I wrote a detailed paper on this subject. Gradually, my interest to the 
problem of best approximation turned into the interest to the methods
that allowed to solve it. And one of them was exactly the method of linear programming.
G.~P.~Akilov advised me to discuss these questions with 
G.~Sh.~Rubinshtein. During our conversations,
G.Sh.\ complemented L.V.'s talks with accounts of 
the related work of  other
mathematicians. Undoubtedly, at that time 
G.Sh.\ was one of the best experts
in linear programming and this circle of L.V.'s ideas in general.
It was not until somewhat later that we learned about the American works
(the simplex method). For us, the principal method was the method of
``resolving multipliers.'' It fitted, as a special case,
in what we called the simplex method, 
but our understanding of the term was wider 
than the American one; the classical Dantzig's simplex method was also
a special case of this, more general, class of methods. Unfortunately,
as it is often the case, the Russian terminology was not well
thought-out and fixed, and the term ``simplex method'' admits 
a lot of various interpretations.

The school of numerical methods of linear programming in the USSR
was extremely strong, and it is undoubtedly the  merit of
L.V.\ and two his principal assistants of the first generation,
V.~A.~Zalgaller and G.~Sh.~Rubinshtein. Later, an important
contribution was made by I.~V.~Romanovsky
and his group, V.~L.~Bulavsky,
D.~B.~Yudin and 
E.~G.~Golshtein in Moscow. Afterwards, with the development of 
computer techniques, numerical solution of problems of 
any reasonable dimension became available.

\subsection{Kantorovich metric}

Once, in the spring of 1957, G.~Sh.~Rubinshtein told me that
he had finally understood how one can use L.V.'s theorem on the
Monge problem (now it is called the Monge--Kantorovich problem)
proved in his 1942 note in {\it Doklady}; namely, how one can
use the Kantorovich metric, i.e., the optimal value of the objective
functional in the mass transportation problem, for introducing a norm
in the space of measures, and how L.V.'s optimality criterion becomes a 
theorem on the duality of the space of measures with the
Kantorovich metric and
the space of Lipschitz functions. In fact, this was
an important methodological remark, since the metric itself had been
already described in L.V.'s note. But it is the paper by L.V.\ and
G.Sh., appeared in {\it Vestnik Leningrad Univ.} in 1958, in the volume dedicated
to G.~M.~Fikhtengolts, that contained the general theory of this,
now famous, metric, 
which is sometimes called the Kantorovich--Rubinshtein
metric, or the transportation metric. By the way, in the same volume
I published my first paper, written jointly with my first scientific advisor
G.~P.~Akilov, devoted to a new definition of Schwartz distributions;
and  we also considered this metric, which had just appeared, 
as one of the  examples. 
A less frequently remembered fact is that the
same paper by L.V.\ and G.Sh.\ 
contained a criterion of optimality of a transportation plan formulated in 
the dual terms
of Lipschitz functions or potentials. Since then I turned into a
permanent propagandist of this remarkable metric, and convinced very
many mathematicians, here and abroad, of the priority of L.V.\ and
the importance of this metric.
It has been rediscovered a lot of times 
(e.g., by L.~N.~Wasserstein or 
D.~Ornstein, who did not know
about L.V.'s work) and thus has many names; 
the method of introducing this metric is known, e.g.,
as coupling, the method of fixed marginal measures, etc. It
has wide applications
in mathematics, especially in ergodic theory, mathematical statistics,
operator theory, and, in the last years, the theory of differential
equations, geometry, and, of course,
statistical physics (see the 2007 Addendum below).
Books are written about it, which
still do not exhaust all its aspects. The Lev\`y--Prokhorov--Skorokhod
metric, well-known in probability theory, is its close relative. 
The possibility of its further generalization to a wide
range of optimization problems was understood somewhat later; this is 
the subject of my 1970 paper in {\it Uspekhi}, which is further developed
in a joint paper with M.~M.~Rubinov.

Simultaneously, in 1970, I applied this metric in an important problem
of measure theory and ergodic theory (in the theory of decreasing
sequences of measurable partitions). There I had to consider
the (wild at first sight) infinite iteration of this metric
(a ``tower of measures''). At about the same time, D.~Ornstein
rediscovered this metric and  introduced it into ergodic theory under
another name (``$\bar d$-metric'').

The history of this metric and all related topics is an excellent
example of how an applied problem (in this case, the mass transportation problem)
can lead to introducing a very useful purely mathematical notion.

\bigskip\noindent
{\bf 2007 Addendum.} The last five or six years have seen
a qualitative jump in the development and applications of the
mathematical theory of mass transportation problems, which we cannot but
mention here. Instead of a rather modest progress in the theory
of Kantorovich's transportation metric in the 70s--80s and the appearance of
several papers a year on this subject (not including many papers of purely
applied or computational nature), several hundred
papers and a number of books were published starting from the late 90s.
As an example we can mention the recent (December 2006) six hundred page survey book 
{\it Optimal Transport, Old and New} by
C.~Villani with six hundred references, mainly to papers of the last
years. There are several reasons for this outburst. The first one is the
sharp increase in the area of application of mass transportation methods:
now it includes not only the traditional fields 
such as mathematical economics, statistics,
probability theory, ergodic theory, but also partial differential equations,
differential geometry, hydrodynamics, some areas of 
theoretical physics, etc. The second one is that the notion of the Kantorovich
metric itself has been generalized: its $p$-analog 
$k_p$, $p\geq 1$, $p \in \mathbb R_+$ (which coincides with the ordinary Kantorovich
metric for $k=1$), which did not attract much attention until recently,
has come into use. The quadratic transportation metric $k_2$
has proved to be the most important; the simplex of probability measures
equipped with this metric apparently can become one of the most
important and useful infinite-dimensional manifolds.
Fruitful connections with the Monge--Amp\`ere problem,
geometric measure theory, Ricci flows (which have become very popular), and other fields
made the study of mass transportation problems one of the central areas
of the modern analysis. Another, quite recent, facet of this subject
concerns already the Kantorovich--Rubinshtein norm.
Note that the $p$-metrics for
$p>1$ do no longer generate any norm on the space of measures. However,
quite recently it became known that the Kantorovich--Rubinshtein norm 
on the space of measures of an arbitrary 
(not necessarily compact) metric space has a simple characteristic property: 
it turns out that it
is exactly the {\it maximal norm} on the vector space of real 
compactly supported measures with bounded variation
on a separable metric space,
in the class of all norms that agree with the metric;
the latter condition means that the norm of the difference of two delta measures
is equal to the distance between the corresponding points:
$\|\delta_x - \delta_y\|=\rho(x,y)$. This provides a
clear geometric interpretation of the norm and a link to root
polytopes of Lie groups and so-called rigid metric spaces, introduced
and investigated in the recent paper
J.~Melleray, F.~V.~Petrov, and A.~M.~Vershik,
Linearly rigid metric spaces and Kantorovich type norms,
{\it C. R. Acad. Sci.} {\bf 344}, no. 4, p. 235 (2007).

I would like to turn once again to the history of the question, being
forced to do this by the striking unanimity in using incorrect terminology
that has gradually seized the whole literature, especially in the West.
On the one hand, L.V.'s priority in the mathematical formulation
of the mass transportation problem, the definition of the transportation
metric, and the introduction of the optimality criterion was universally recognized long 
ago. Usually, this mass transportation problem of Kantorovich
is called the Monge--Kantorovich problem, because G.~Monge, indeed,
considered the plane problem of 
transporting a pile of soil to an excavation fill
as a transportation problem with a similar estimation of costs.
By the way, L.V.\ did not learn about this work till the mid 50s,
when Monge's collected works were published on the occasion
of his 200th birthday; 
L.V.\ wrote a small note in {\it Uspekhi} in which he explained how to
fit Monge's problem into the required framework and how the
optimality criterion works in this case. Monge's voluminous work contained no
mention of the metric and, all the more so, the optimality criterion, 
so that the term ``Kantorovich metric'' as applied to the transportation
metric is absolutely justified, and he had
no other predecessors in considering this problem. Perhaps, 
the name of the mass transportation problem itself, the
Monge--Kantorovich problem, can also be accepted,
though with great reserve,
because in Monge's times there were no general metric spaces,
let alone that the very idea
to set the problem in such a generality and apply the duality method
borrowed from functional analysis is highly nontrivial, and it was an outstanding
achievement of L.V. We should also take into account that for a number
of objective reasons (the War, the separation of the Soviet mathematicians from the
West, etc.) the 1942 note in {\it Doklady}, as well as other L.V.'
works in mathematical economics, for a long time remained unknown to
the mathematical world. As written above, for a long time L.V.\
did not popularize his work in mathematical economics because of the pressure
of the Soviet obscurant ideological censorship.
Because of this forced
delay, or for other reasons, a number of authors working in various fields,
sometimes not connected with mathematical economics, rediscovered and applied this
metric in some or other concrete situation (most often, in
less general setting than that of L.V.), not knowing about the works
by L.V.\ and his successors. Here one can mention L.~Wasserstein,
D.~Ornstein, and others. Sometimes, it was not an easy task to convince
such an author or his colleagues that his discovery was already known,
but eventually one managed to do this. But it is certainly
incorrect to attribute this metric, e.g., to 
L.~Wasserstein or to other rediscoverers, as it is being done
by many authors. Moreover, this is also incorrect
from the formal point of view: none of the subsequent works
that contained, explicitly or quite implicitly, the definition of this metric,
was of such a generality, clearness, and fundamentality that made this
1942 paper of L.V.\ classical.

To repair this injustice
is difficult but necessary.

In 2004, an international conference dedicated to the 90th
anniversary of L.V.\ took place in St.~Petersburg. 
Its program and many of the talks are published in {\it Zapiski
Nauchnyh Seminarov POMI}, vol.~312 (2004); the English translation
of this volume appeared in {\it Journal of Mathematical Sciences (New York)},
vol.~133, no.~4 (2006).
It also contains a reproduction of L.V.'s 1942 note and a number of comments
on this circle of problems.

\subsection{Connections with calculus of variations and Lagrange multipliers}

Linear and convex programming is a natural generalization of the theory of Lagrange multipliers
to nonregular problems (problems on polyhedral domains, or, as we would
say now, on manifolds with corners). The fact that resolving multipliers
are a generalization of Lagrange multipliers was noticed by L.V.\
at the very beginning. Nonclassical multipliers appeared also in
other areas, first of all, in the theory of optimal control in
L.~S.~Pontryagin's school. This theory also generalized conditional
variational problems to the case of nonregular constraints, and thus
it should be compared with problems of (in general, nonconvex, but
in most important cases, convex) infinite-dimensional programming.
This link did not become clear till somewhat later. 
One should say that aesthetically Pontryagin's theory is inferior
to L.V.'s one, though the former is substantially more complicated
(only because its problems are originally infinite-dimensional).
There is a lot of literature on the connections of linear and
convex programming with optimal control. However, for a number of reasons,
these connections have not been elaborated to a sufficiently deep level.
The principal reason is that the  form in which 
one usually considers problems of optimal control is
insufficiently invariant. An intermediate
position between classical calculus of variations and optimal control,
lying closer to geometry and the theory of Lie algebras, is occupied by
nonholonomic problems. They also involve nonclassical constraints,
like convex programming and optimal control, but 
of another, smooth, kind.
I started to study nonholonomic problems in the mid 60s, when I begun to think about the
works on invariant formulations of mechanics
(by V.~I.~Arnold, C.~Godbillon, J.~E.~Marsden, etc.),
much popular at that time. Having seen in nonholonomic mechanics, which is
a stepdaughter of classical mechanics, a nontrivial optimization
problem, I understood how to state it in modern form. In those years, at LOMI
we had an educational seminar for young people
(its participants included L.~D.~Faddeev, B.~B.~Venkov, me, and others), at which we studied
differential geometry, representation theory, Lie groups, and all other
things. One day it emerged
by pure accident that L.D.\ also thought about nonholonomic mechanics,
and we decided to get to the bottom of this matter together. We wrote first
a brief note for {\it Doklady} and then a large paper on
the invariant form of Lagrangian and, in particular, nonholonomic mechanics.
These papers are still abundantly cited; they contain a dictionary
of correspondences between the terms of differential geometry and
the notions of classical mechanics. Now this subject has become very popular;
it is a remarkable intermediate between classical and nonclassical 
calculus of variations. In this setting, Lagrange multipliers appear in yet another new form, namely,
as variables corresponding to constraints and consequences
(Lie brackets) of all orders. Here one also cannot help remembering
L.V.'s resolving multipliers.

\subsection{Linear models and Markov processes}

Since in the 60s L.V.\ worked intensively on economic models, not
necessarily related to optimization, one cannot but mention,
at least briefly, the connections of the theory of models of
economic dynamics (works by J.~von Neumann, W.~Leontief, L.V., and others)
with dynamical systems. Here I would only like to emphasize one
connection, still insufficiently studied; namely, that these linear economic models
are directly related to a particular type of Markov processes, in which
a special role is played by the notion of positivity in the set of states. 
Turnpike-type
theorems and Markov decision-making processes, as well as
theories
of multivalued maps, continuous choice problems, etc.,
are most directly related to
this circle of problems.  Apparently,
these problems are now losing their importance for applications, but
they are undoubtedly interesting from the viewpoint of mathematics,
like any theories of multivalued and positive maps. Recall that
before the war L.V.\ created the theory of partially ordered spaces (K-spaces),
which soon closed on itself and ceased to interest him and anyone
who did not work directly in this field. But partial ordering
in a more
general sense has always been a subject of special interest for
mathematicians of the Leningrad and Ukrainian schools.

\subsection{Globalization of linear programming}

Attracting ideas from topology and differential geometry 
led also to another synthesis, namely, to the notions of fields of polytopes, cones,
etc., playing an important role in optimal control, Pareto optimum
(S.~Smale's conjecture and works of Y.-H.~Wan and
A.~M.~Vershik--A.~G.~Chernyakov), etc. I mean problems with a smooth parameter
that ranges over a manifold at each point of which there is a problem
of linear programming. Fields of polytopes, or fields of problems,
also arise in the theory of smooth dynamical systems. Another area
in which one uses similar tools but with another goal is 
the problem of estimating the number
of steps in various versions of the simplex method (works of S.~Smale, 
A.~M.~Vershik--P.~V.~Sporyshev,
etc.); solving this problem involved
ideas of integral geometry (``Grassmanian approach'').
The obtained estimates gave another evidence of the practical value
of the simplex method and the method of resolving multipliers. 
In the 80s, a strong
impression was produced by the works of L.~G.~Khachiyan and N.~Karmarkar
which gave a polynomial (in a sense) uniform (in the class of problems)
bound on the complexity of the ellipsoid method for solving
problems of linear programming. Nevertheless, this method in
no respect replaced various versions of the simplex method. The bounds mentioned
above lead to linear or quadratic bounds on the complexity that hold only
statistically. On the whole, it is still (2001) unknown whether
the problem of linear programming 
belongs to the class P of polynomial problems.

\subsection{Linear programming and computational methods}

Another research direction that was initiated by L.V.\ and has not been properly developed
is the problem of applying linear programming to
approximate solution of problems
of mathematical physics (two-sided bounds on linear functionals of
solutions). L.V.'s paper on this subject (1962)
contained a very fruitful idea,
and several related works were carried out at the Leningrad University.
L.V.'s approach can also be considered as an alternative approach to 
noncorrect problems. This problem is very important in geophysics,
and L.V.\ discussed it with V.~I.~Keilis-Borok.

\section{L.V.\ and personnel training}

One of the important initiatives of L.V.\ of that period 
(50s--60s) was to start
training of specialists in mathematical economics. 
Already in the 50s, L.V.\ had a number of students 
and pupils working in this field, but in comparison with
other his numerous studies, the number of pupils in this field
was rather small. In earnest this
training began in 1959, when the so-called 
sixth year of education was organized at the Department of Economics
of the Leningrad State University; it allowed students to 
become acquainted with mathematical
economics and L.V.'s ideas. Many well-known economists completed
this sixth year: A.~A.~Anchishkin, S.~S.~Shatalin, I.~M.~Syroezhin, etc.
This course (which existed only one year) became the center of mathematical 
retraining of economists. It is not out of place to recall that
the most part of prominent economists of the 70s--90s have
in some way learned from  L.V.\ or communicated with him. 
Among those
closest to him, I would like to mention
A.~G.~Aganbegyan and V.~L.~Makarov. Soon, in 1959, the Chair
of Economic Cybernetics was organized at the Department of Economics. 
At the first stage, a very active role in the organization of the new chair
was played by V.~V.~Novozhilov, an old companion of L.V.\ in
battles with conservative economists and the author of fascinating economic conceptions.
The mathematicians that most
actively participated in the organization and teaching in the first years
are V.~A.~Zalgaller and, somewhat later, L.~M.~Abramov. Among 
political economists one should mention  I.~V.~Kotov, 
the future Head of the Chair, and V.~A.~Vorotilov, the then Head of the Department
of Economics, as well as I.~M.~Syroezhin, the head of the Laboratory
of Mathematical Economics.
As a matter of fact, the mathematical ``intrusion''
into the Department of Economics
had far-reaching consequences not only for economic
cybernetics, but 
for the Department as a whole. Since then mathematics held a firm place
at the Department of Economics, 
and mathematical education there became comparatively good;
mathematical courses were mainly taught by professors from the Department of Mathematics
and Mechanics and at the same level. L.V.'s flying visits from Novosibirsk,
though not very frequent, were nevertheless very fruitful: the most important
decisions concerning the new speciality were, to a certain extent, taken
on his behalf. Somewhat later (already after L.V.'s departure for Novosibirsk,
but with his participation), 
the same was done at the Department of Mathematics
and Mechanics: first (in 1961--1962), the speciality ``operation research'' was
organized in the bosom of the Chair of Computational Mathematics,
and then, in 1970, the Chair of Operation Research was created.
The principal role in the development of this chair was played by M.~K.~Gavurin and
I.~V.~Romanovsky, who, since the 60s, taught his seminar on optimization
with a bias to computational problems.

Economic cybernetics quickly found its niche. The necessity to
mathematize and renew the dilapidated (of course, this fact
was not officially recognized) economic science, to study the functioning
of economic structures and to optimize them quite naturally required 
a new type of specialists. The new
chairs of economic departments were intended to train such specialists.

At the same time, strange as it may seem, there were certain difficulties
in finding the place of this speciality 
within mathematics itself.  The new speciality organized 
at the Department of Mathematics and Mechanics of the Leningrad
State University was
one of the first in the country, it was organized
almost simultaneously with a similar speciality at the Novosibirsk University.
The difficulties were
in the fact that, important as they were, the
models and methods of mathematical economics did not constitute a new area
of theoretical mathematics. The mathematical aspects of the theory
created by L.V., or W.~Leontief, or J.~von Neumann, fitted well,
on the one hand, in the framework of functional (more exactly, convex)
analysis, the theory of inequalities, etc., and, on the other hand 
(from the practical point of view), in the framework of
the theory of computational methods
(in which L.V.\ also was one of the leaders) for solving extremal problems.
As to the theory of linear programming, it was a spectacular and natural
generalization of classical methods (Lagrange multipliers,
conjugate problems, 
duality, etc.). Anyhow, all this (plus optimal control)
could be called new directions, new areas, but not a new mathematical
science, as was the case with economic cybernetics, or, more exactly,
mathematical economics, in the framework of economic science.
As was said above, at first (since 1962) the speciality ``operation research''
was taught at the Chair of Computational Mathematics.
I remember well one of the conversations of L.V.\ with the then Head
of the Department, for which I was invited (being yet a postgraduate
student). The Head, who did not thoroughly realize the purely
mathematical weight of the new area, tried to persuade me to work
entirely on the mathematical questions related to L.V.'s ideas,
and L.V., who supported my candidacy for a position at the chair, answered that 
it was not quite enough for me from the point of view of ``pure mathematics.''
After much foot-dragging, mainly of nonscientific nature, I was nevertheless
hired by the Department, but not by the Chair of Mathematical
Analysis, from which I had graduated and obtained my Ph.D. degree, but
by the Chair of Computational Mathematics, specially for teaching
the new speciality. There was indeed some vagueness
in the position of the Chair and the speciality itself, since
it had not its sharply defined specific characteristics (say, as the Chair of Algebra,
or Geometry, or even Computational Mathematics) and thus was compelled to
become interdisciplinary and partly applied. Its research and teaching area had
intersections with the areas of various chairs: the Chair of Differential
Equations (variational problems), the Chair of Mathematical
Analysis (convex and functional analysis), the Chair of Algebra
(discrete mathematics), and, of course, the Chairs of Computational
Mathematics and Software Engineering). But its own area was not sufficiently wide
to become an independent theoretical mathematical specialization. 
This predetermined both strong and weak sides of the future chair
and speciality. Let me add, in parentheses, that I was, and remain,
against partitioning mathematical departments into chairs; this old German
tradition did not survive in any of the leading mathematical countries.
Now (and for a long time) it only impedes necessary changes in the system
of mathematical education. As far as I know, 
the efficiency of education at the Department of Mathematics and Mechanics
has not been seriously studied,
but I am afraid that the form of education that has not undergone any
changes for so long just cannot turn out to be good. Because of this,
the speciality and the Chair did not attract particularly strong
students.

The situation in theoretical economics was quite different. There, the new
ideas had attracted the most fresh and healthy forces, and subsequently
L.V.\ became the doubtless leader and teacher of a whole pleiad
of our economists. It would not be an exaggeration to say 
that all a bit educated modern Russian economists are either L.V.'s pupils, or pupils of his pupils,
or have somehow absorbed his ideas.
Of course, this is an
important subject
for a special historical investigation. It is difficult
for me to tell about the Novosibirsk and Moscow periods of his
pedagogical and scientific activity: this is quite another epoch
(and even two epochs), which is apparently not similar to the Leningrad period. 

\section{Several personal reminiscences}

The personality of L.V., his qualities as a teacher and a scientist,
are worth a separate paper. Here I will restrict myself to
several remarks.

\medskip
{\bf 1.} My first meetings, conversations, and contacts with him 
especially impressed me, and my friends, with the speed with which L.V.\
grasped what was being said to him, 
anticipating and immediately calculating what was going to be said.
Later I read the same about J.~von Neumann;
by the way, before the war he had a correspondence with L.V.\
on a circle of problems related to partially ordered spaces. The very first
L.V.'s works (joint with E.~M.~Livenson) on descriptive set theory, from which
his fame had begun, impressed the Moscow mathematicians, who had been
working in this field for a long time, with technical skills and
deep insight. His versatility and the exact understanding of the
essence, whatever was the matter, were also striking. The speed and
depth of his mathematical thought were on the edge of human capabilities 
(at least, those known to me).

I remember a discussion of a series of papers of American authors 
on automata theory, very popular at the time, that took place
in the 60s
at a seminar at the Leningrad House of Scientists. In particular, L.V.\ commented
upon the paper ``Amplifier of intelligence''
by W.~R.~Ashby,
in which the author substantiated the obvious idea of
necessity to speed up the mental process. L.V.: ``Of course,
the speed of thought is different for different people, but
it can differ from the average by three, well five, but not by
1000 times.'' I think L.V.'s coefficient was much greater than five.

\medskip
{\bf 2}. At the same time, he read his lectures at a slow but quite uneven pace, 
reacting very vividly to questions. Each lecture began with the
sacramental question ``Are there any questions on the previous lecture?'',
which used to be pronounced by a rolling loud voice. But sometimes,
during a lecture, this voice lowered almost to a whisper. 
At seminars, he frequently slept, but nevertheless, by some miracle, 
he interrupted the speaker at appropriate moments, looking far ahead
of what had been said.  His comments were always useful and instructive.

\medskip
{\bf 3.} But as to lectures of fundamental importance, he 
delivered them  brilliantly.
He was an extremely experienced polemicist and could always find exact
objections that were at the heart of the matter. I remember well a number of his
speeches mentioned above. It's a pity that at that time
there was no videotape recording.

\medskip
{\bf 4.} According to my observations, his attitude towards mathematics varied.
Before the war and in the first post-war years, he undoubtedly belonged
to the small number of leaders of functional analysis (other ones were
I.~M.~Gelfand, M.~G.~Krein). This became especially clear after his
famous paper ``Functional analysis and applied mathematics''
published in {\it Uspekhi}, for which he was awarded the Stalin Prize, 
very important for his subsequent stability in
troublesome times. His well-known book 
written jointly with G.~P.~Akilov summarized the
activity of the Leningrad school of functional analysis.
Later, having turned to economics, he went slightly away from
mathematics, but, in my opinion, he well understood that the construction of the
functional analysis that had been created and successfully developed
in Leningrad in the 30s--50s was already almost completed, 
and tried to popularize in Leningrad new research directions.
I remember well his interest to the theory of Schwartz distributions;
once, in 1956, on his and G.~P.~Akilov's request, I gave a talk
at the Fikhtengolts--Kantorovich seminar on various definitions
of generalized functions, and one of the first definitions was that
of L.V.\ from a note published in {\it Doklady} in 1934, earlier than
Sobolev's works, etc.! Later, on numerous occasions he told me about the role
of I.~M.~Gelfand in mathematics and expressed regret that the latter
was not yet elected to the Academy.

It seemed L.V.\ was sorry that after the 50s he had in fact left 
mathematics, but, in my opinion, his choice between mathematics and economics
was predetermined. 

\medskip
{\bf 5.} But L.V.\ could also serve as an excellent example of a scientist
who could be called an ``applied mathematician.'' His flair for applied 
problems and the most extensive contacts with engineers, military
specialists, economists made him very popular among those
who applied mathematics. In his own words,  he felt himself
not only a mathematician, but also an engineer. His successful work
in engineering, programming, industrial computations is an excellent illustration
of this thesis.

\medskip
{\bf 6.} In the professional environment he was almost always surrounded
by general admiration and attention. 
If he was in good form, his appearance at seminars, lectures
at once animated the atmosphere, 
``brownized'' it, as one says. I believe that this was acknowledged by everybody ---
by his friends as well as by his enemies. In the last years, having already 
gone away from mathematics, he was friends with the leading Moscow
mathematicians of the next generation: V.~I.~Arnold, S.~P.~Novikov, etc.
I hope that they will sometimes write about their conversations with him.
\medskip

In conclusion of this essay, I wish to say that we (my generation
of mathematicians grown up in Leningrad), and me personally, were
incredibly lucky with our teachers and we were lucky to witness
and even slightly participate in creating new scientific directions,
and to learn from their creators. Here I would like to emphasize the role of
L.V., which is not yet thoroughly understood and appreciated.
At first sight, as he himself used to say (but here one should make
a natural allowance for
internal and external censorship), 
his theories were adjusted to planned economy. But  
this is only an outward appearance. The main things, such as 
taking into account hidden parameters (rent), a unified approach to
constraints (labor is only one of them), and all the ensuing 
consequences, make his economic applications universal
and indispensable now. Upon the whole, the main result
of the great experiment of Kantorovich is that
the most modern mathematical tools can not only be 
successfully applied in economics, but they can serve as a basis for 
creating new purely economical conceptions.
This does not mean that his conclusions
will fully work today. But this undoubtedly means, and in this respect
L.V.\ was perhaps the first (J.~von Neumann did not work in economics at such a
deep level as L.V.), that the mathematician's talent 
can radically alter and transform economic thought. Most
unfortunately, L.V.\ did not live until the 90s, when his
experience, intuition, and authority could be used with much greater effect
than in Soviet times. I have no doubt that he would be able to caution
the reformers, whose theoretical (and even practical)
skills were at an insufficiently high level (which made them listen to
dubious advises) against serious mistakes. Alas, when it was necessary,
in Russia there turned out to be no experienced economist of
L.V.'s magnitude.

\end{document}